
\documentstyle{amsppt}
\baselineskip18pt
\magnification=\magstep1
\pagewidth{30pc}
\pageheight{45pc}

\hyphenation{co-deter-min-ant co-deter-min-ants pa-ra-met-rised
pre-print pro-pa-gat-ing pro-pa-gate
fel-low-ship Cox-et-er dis-trib-ut-ive}
\def\leaderfill{\leaders\hbox to 1em{\hss.\hss}\hfill}

\

\def\ldescent#1{{\Cal L (#1)}}
\def\rdescent#1{{\Cal R (#1)}}

\def\ti{\widetilde}

\def\b0{\text{\bf 0}}

\def\lan{{\langle}}
\def\ran{{\rangle}}

\def\zed{{\Bbb Z}}

\def\boxit#1{\vbox{\hrule\hbox{\vrule \kern3pt
\vbox{\kern3pt\hbox{#1}\kern3pt}\kern3pt\vrule}\hrule}}
\def\rabbit{\vbox{\hbox{\kern0pt
\vbox{\kern0pt{\hbox{---}}\kern3.5pt}}}}

\def\tableau#1{
        \hbox {
                \hskip -10pt plus0pt minus0pt
                \raise\baselineskip\hbox{
                \offinterlineskip
                \hbox{#1}}
                \hskip0.25em
        }
}

\def\tabCol#1{
\hbox{\vtop{\hrule
\halign{\strut\vrule\hskip0.5em##\hskip0.5em\hfill\vrule\cr\lower0pt
\hbox\bgroup$#1$\egroup \cr}
\hrule
} } \hskip -10.5pt plus0pt minus0pt}

\def\CR{
        $\egroup\cr
        \noalign{\hrule}
        \lower0pt\hbox\bgroup$
}



\def\blank#1#2{
\hbox to #1{\hfill \vbox to #2{\vfill}}
}


\def\strut{\vrule height10pt depth5pt width0pt}

\topmatter
\title Leading coefficients of Kazhdan--Lusztig polynomials and fully
commutative elements
\endtitle

\rightheadtext{Leading coefficients of Kazhdan--Lusztig polynomials}

\author R.M. Green \endauthor
\affil Department of Mathematics \\ University of Colorado \\
Campus Box 395 \\ Boulder, CO  80309-0395 \\ USA \\ {\it  E-mail:}
rmg\@euclid.colorado.edu \\
\newline
\endaffil


\abstract
Let $W$ be a Coxeter group of type $\widetilde{A}_{n-1}$.  We show that 
the leading coefficient, $\mu(x, w)$, of the Kazhdan--Lusztig polynomial 
$P_{x, w}$ is always equal to $0$ or $1$ if $x$ is fully commutative (and 
$w$ is arbitrary).
\endabstract

\endtopmatter

\centerline{\bf Preliminary version, draft 3}

\head Introduction \endhead

In their famous paper \cite{{\bf 4}}, Kazhdan and Lusztig showed how to 
associate to an arbitrary Coxeter group $W$ a remarkable family of
polynomials, $\{P_{x, w}(q) : x, w \in W\}$, which are now known as
Kazhdan--Lusztig polynomials.  The Kazhdan--Lusztig polynomials are of 
key importance in algebra and geometry.  For example, they are intimately 
related to the geometry of Schubert varieties, and they are necessary 
for the statement of Lusztig's famous conjecture \cite{{\bf 7}} regarding the 
characters of irreducible modules of reductive algebraic groups in 
characteristic $p > 0$. 

The polynomial $P_{x, w}$ is zero unless $x \leq w$ in the Bruhat order on $W$.
If we have $x < w$, then 
$P_{x, w}(q)$ is of degree at most $(\ell(w) - \ell(x) - 1)/2$, where $\ell$
is the length function on the Coxeter group.  The cases where this degree
bound is achieved are of particular importance, and in such cases, the
leading coefficient of $P_{x, w}(q)$ is denoted by $\mu(x, w)$.  The 
$P_{x, w}(q)$ and $\mu(x, w)$ are defined by recurrence relations and are
very difficult to compute efficiently, even for some moderately small groups.

As well as playing a key role in the computation of the polynomials
$P_{x, w}(q)$, the $\mu$-function is also an interesting object in its own
right.  Empirical evidence had suggested that in the special
case of the symmetric group (Coxeter type $A$), the values $\mu(x, w)$
were always equal to $0$ or $1$.  This hypothesis, which was known as
the {\it $0$--$1$ Conjecture}, was proved to be false by McLarnan and
Warrington \cite{{\bf 10}}.  However, as pointed out by B.C. Jones \cite{{\bf 3}, \S1},
the values $\mu(x, w)$ in the symmetric group case are indeed equal to 
$0$ or $1$ if one restricts $x$ and/or $w$ to certain (interesting) 
classes of permutations.  For example, we have $\mu(x, w) \in \{0, 1\}$
if one of the following holds:
\item{(a)}{the symmetric group $S_n$ in question satisfies $n \leq 9$ 
\cite{{\bf 10}};}
\item{(b)}{$w$ is a Grassmannian permutation, meaning that there is at most 
one decreasing consecutive pair of entries in the one-line notation for $w$
\cite{{\bf 6}};}
\item{(c)}{$a(x) < a(w)$ \cite{{\bf 13}}, where $a$ is Lusztig's $a$-function 
from \cite{{\bf 9}};}
\item{(d)}{$w$ corresponds to a smooth Schubert variety \cite{{\bf 5}}.}

\noindent (In case (d), the polynomials $P_{x, w}$ themselves are equal to 
$0$ or $1$.)  The main result of \cite{{\bf 3}} is that if $W$ is a finite Weyl
group, then $\mu(x, w)$ is always $0$ or $1$ if $w$ is a Deodhar element;
the latter are a subclass of the fully commutative elements.

In this paper, we make a restriction on $x$ rather than $w$, and we work in 
the larger Coxeter group $W = W(\ti{A}_{n-1})$ of type $\ti{A}_{n-1}$.  
This group, which is infinite, naturally contains the symmetric group 
$S_n$ as a subgroup, and it can be thought of as a certain group of 
periodic permutations of the integers \cite{{\bf 8}}.  We will restrict $x$
to be a fully commutative element in the sense of Stembridge \cite{{\bf 12}}.
There are infinitely many such elements in the group $W$, and for each choice
of $x$, there are infinitely many $w$ with $x \leq w$.  Our main result is
that if $x, w \in W$ and $x$ is fully commutative, then 
$\mu(x, w) \in \{0, 1\}$.

Computational evidence suggests that there may be a large class of Coxeter
groups (possibly all of them) for which $\mu(x, w) \in \{0, 1\}$ when $x$
is fully commutative.  We intend to return to this question elsewhere.

\head 1. Definitions \endhead

Let $n \geq 3$, and let $W = W(\ti{A}_{n-1})$ be Coxeter group of type 
affine $A_{n-1}$ with the distinguished set of generating involutions $$
S = \{s_i : 0 \leq i < n \}.$$  In other words, $W$ is 
given by the presentation $$
W = \lan S \ | \ (s_i s_j)^{m(i, j)} = 1 \text{ for } m(i, j) < \infty \ran
,$$ where we define $m(i, i) = 1$, $m(i, j) = 3$ if $|i - j| = \pm 1 \mod n$, 
and $m(i, j) = 2$ otherwise.
The elements of $S$ are distinct as group elements, and $m(i, j)$ 
is the order of $s_i s_j$.  The subgroup of $W$ generated by $S \backslash
\{s_0\}$ is isomorphic to the symmetric group $S_n$; one isomorphism is
given by the correspondence $s_i \leftrightarrow (i, i+1)$.

Each element $w \in W$ can be written as a word $w = s_{i_1} s_{i_2} 
\cdots s_{i_k}$ in the generators $S$.  If this $k$ is minimal for a given $w$,
then we call $k$ the {\it length} of $w$ and write $k = \ell(w)$.  We call 
$s_{i_1} s_{i_2} 
\cdots s_{i_k}$ a {\it reduced expression} for $w$.  More generally, we call
a product $w_1w_2\cdots w_k$ of elements $w_i\in W$ 
{\it reduced} if \newline $\ell(w_1w_2\cdots w_k)=\sum_i\ell(w_i)$.
We write $$
\ldescent{w} = \{s \in S : \ell(sw) < \ell(w)\}
$$ and $$
\rdescent{w} = \{s \in S : \ell(ws) < \ell(w)\}
.$$  The set $\ldescent{w}$ (respectively, $\rdescent{w}$) is called the 
{\it left} (respectively, {\it right}) {\it descent set} of $w$.
We call a (left or right) descent set {\it commutative} if it consists of
mutually commuting generators.

There is a natural partial order, $\leq$, on $W$ called the {\it (strong) 
Bruhat order} and denoted by $\leq$.  It is characterized 
by the property that if $t_1 t_2 \cdots t_r$ is a reduced expression for 
$w \in W$, then the set of
elements $x \leq w$ are precisely those of the form $
x = t_{i_1} t_{i_2} \cdots t_{i_k}
,$ where we have $1 \leq i_1 < i_2 < \cdots < i_k \leq r$.

For each $x, w \in W$, Kazhdan and Lusztig \cite{{\bf 4}} defined a 
polynomial $P_{x, w} \in \zed[q]$.  We have $P_{w, w} = 1$, and 
unless $x \leq w$ in the Bruhat order, we have $P_{x, w} = 0$.  The other
possibility is that $x < w$, in which case
$P_{x, w}$ is a polynomial in $q$ of degree at most 
$(\ell(w) - \ell(x) - 1)/2$.  We define $\mu(x, w)$ to be the (integer) 
coefficient of $q^{(\ell(w) - \ell(x) - 1)/2}$ in $P_{x, w}$; this can only 
be nonzero if (a) $x < w$ and (b) $\ell(x)$ and $\ell(w)$ 
have opposite parities.  We write $x \prec w$ to mean that both 
$\mu(x, w) \ne 0$ and $x < w$.  Using these facts, one can define the
Kazhdan--Lusztig polynomials recursively via a recurrence formula proved
in \cite{{\bf 4}}.  For this, we may assume that $w \ne 1$, for otherwise
we have $x = 1$ and $P(x, w) = 1$.  We choose $s$ such that 
$\ell(sw) < \ell(w)$ (take
$s$ to be the first letter in a reduced expression for $w$) and we define
$v = sw$.  We may then compute $P_{x, w}$ recursively by the formula $$
P_{x, w} = q^{1-c} P_{sx, v} + q^c P_{x, v} - 
\sum_{{z \prec v} \atop {sz < z}} \mu(z, v) \, 
q^{(\ell(w) - \ell(z))/2} P_{x, z}
,$$ where we define $c = 0$ if $x < sx$ and $c = 1$ otherwise.
Notice that a knowledge of the values of the $\mu$ function is important
for computing the polynomials.  However, as we shall see, it is possible to
compute $\mu$ values without first computing the corresponding polynomials.

We call an element $w \in W$ {\it complex} if it can be written 
as a reduced product $x_1 w_{ij} x_2$, where $x_1, x_2 \in W$, $|i - j| 
= \pm 1 \mod n$ and $w_{ij} = s_i s_j s_i = s_j s_i s_j$.
Denote by $W_c$ the set of all elements of $W$
that are not complex.  The elements of $W_c$ are the 
{\it fully commutative} elements of \cite{{\bf 12}}; they are characterized by 
the property that any two of their reduced expressions may be obtained 
from each other by repeated commutation of adjacent generators.

A key concept for this paper is that of a star operation, which was
introduced in \cite{{\bf 4}, \S4.1}.
Let $I = \{s, t\} \subseteq S$ be a pair of noncommuting generators of $W$.
Let $W^I$ denote the set of all $w \in W$ satisfying 
$\ldescent{w} \cap I = \emptyset$.  Standard properties of Coxeter groups 
\cite{{\bf 2}, \S5.12} show that any element $w \in W$ may be uniquely
written as $w = w_I w^I$ reduced, where $w_I \in W_I = \lan s, t \ran$
and $w^I \in W^I$.  There are four possibilities for elements $w \in W$:
\item{(i)}{$w$ is the shortest element in the coset $W_I w$, so $w_I = 1$ and 
$w \in W^I$;}
\item{(ii)}{$w$ is the longest element in the coset $W_I w$, so $w_I$ is 
the longest element of $W_I$ (which relies on $W_I$ being finite);}
\item{(iii)}{$w$ is one of the two elements $sw^I$ or $tsw^I$;}
\item{(iv)}{$w$ is one of the two elements $tw^I$ or $stw^I$.}

The sequences appearing in (iii) and (iv) are called {\it (left) 
$\{s, t\}$-strings}, or {\it strings} if the context is clear.  
If $w$ is an element of an $\{s, t\}$-string, $S_w$, then we define
$^*w$ to be the other element of $S_w$.  If $w$ is not an element of an
$\{s, t\}$-string (in other words, case (i) or (ii) applies) then $^*w$
is undefined.

There are also obvious right handed analogues to the above concepts,
so the symbol $w^*$ may be used with the analogous meaning (with respect
to pair $I$ of noncommuting generators).

\head 2. Preparatory results \endhead

In \S2, we collect some results from the literature for use in the proof
of our main result.

The following key properties of the $\mu$-function 
were proved in \cite{{\bf 4}}.

\proclaim{Proposition 2.1}
\item{\rm (i)}{If $\ell(x) = \ell(w) \mod 2$ or $x \not\leq w$, then we
have $\mu(x, w) = 0$.}
\item{\rm (ii)}{Let $x, w \in W$ be elements of left $\{s, t\}$-strings 
(for the same $s$ and $t$, but possibly different strings).  Then we
have $\mu(x, w) = \mu(^*x, {^*w})$.}
\item{\rm (iii)}{Let $x, w \in W$ be elements of right $\{s, t\}$-strings 
(for the same $s$ and $t$, but possibly different strings).  Then we
have $\mu(x, w) = \mu(x^*, w^*)$.}
\item{\rm (iv)}{If there exists $s \in \ldescent{w} \backslash \ldescent{x}$
then we have either (a) $\mu(x, w) = 0$ or (b) both 
$x = sw$ and $\mu(x, w) = 1$.}
\item{\rm (v)}{If there exists $s \in \rdescent{w} \backslash \rdescent{x}$
then we have either (a) $\mu(x, w) = 0$ or (b) both
$x = ws$ and $\mu(x, w) = 1$.}
\qed\endproclaim


\definition{Definition 2.2}
Suppose that $n$ is even.  Let $$
S_0 = \{s_0, s_2, s_4, \ldots, s_{n-2}\}
$$ and $$
S_1 = \{s_1, s_3, s_5, \ldots, s_{n-1}\}
.$$  We define $I_0$ (respectively, $I_1$) to be the product in 
$W(\ti{A}_{n-1})$ of the generators in $S_0$ (respectively, $S_1$).
(Note that each of $I_0$ and $I_1$ is a product of $n/2$ commuting
generators.)
\enddefinition

\proclaim{Proposition 2.3 (Shi)}
Let $w \in W(\ti{A}_{n-1})$.
If $w \not\in W_c$, then there is a finite sequence $
w = w_0, w_1, \ldots, w_k
$ such that 
\item{\rm (i)}{for each $0 \leq i < k$, we have $w_{i+1} = {^*w}_i$ with
respect to some pair of noncommuting generators $\{s(i), t(i)\}$ depending on
$i$, and}
\item{\rm (ii)}{$\ldescent{w_k}$ is not commutative.}
\endproclaim

\demo{Proof}
This is a special case of \cite{{\bf 11}, Lemma 2.2}.
\qed\enddemo

\proclaim{Proposition 2.4 (Fan--Green)}
Let $w \in W(\ti{A}_{n-1})$.
Suppose that $w \in W_c$.  Then one of the following four
situations must occur:
\item{\rm (i)}{$w$ is a product of commuting generators;}
\item{\rm (ii)}{$n$ is even and $w$ is equal to an alternating product of the
elements $I_0$ and $I_1$ in Definition 2.2;}
\item{\rm (iii)}{we have $w = stv$ reduced for a pair of noncommuting 
generators $I = \{s, t\}$, and we have $tv = {^*w}$ with respect to $I$.}
\item{\rm (iv)}{we have $w = vts$ reduced for a pair of noncommuting 
generators $I = \{s, t\}$, and we have $vt = w^*$ with respect to $I$.}
\endproclaim

\demo{Proof}
This is a restatement of \cite{{\bf 1}, Proposition 3.1.2}.
\qed\enddemo

\proclaim{Lemma 2.5}
Let $w \in W(\ti{A}_{n-1})$ be fully commutative.  
\item{\rm (i)}{The sets $\ldescent{w}$ and $\rdescent{w}$ are commutative 
(in the sense of \S1).}
\item{\rm (ii)}{Suppose $n$ is even, that $\ldescent{w} \in \{S_0, S_1\}$,
and that $\rdescent{w} \in \{S_0, S_1\}$.  Then $w$ is equal to an 
alternating product of $I_0$ and $I_1$.}
\endproclaim

\demo{Proof}
If $s$ and $t$ are noncommuting generators in $\ldescent{w}$, then it follows
using basic properties of Coxeter groups that $w$ has a reduced expression
beginning $sts$, which is incompatible with the hypothesis $w \in W_c$.
A similar argument deals with the case of $\rdescent{w}$, and assertion (i) 
(which is well known) follows.

Assume the hypotheses of (ii).  Suppose in addition that $\ldescent{w} = S_0$,
and assume (for a contradiction) that we are in (iii) of Proposition 2.4.  
In this case, we have $s \in S_0$.  There will be precisely two generators
that do not commute with $t$, both of which lie in $S_0$; let us call
the other such generator $u$.  The hypothesis $\ldescent{w} = S_0$ shows that
$u \in \ldescent{w}$.  Now $w = stv$, being fully commutative, has the 
property that the occurrence of $t$ shown lies to the left of any occurrence 
of $u$, and this property is retained after repeated commutation of adjacent 
generators.  It follows that $w$ has no reduced expression beginning in $u$,
which is incompatible with $u \in \ldescent{w}$ by the Exchange Condition for
Coxeter groups.  This is a contradiction.

Using similar arguments, we see that the case $\ldescent{w} = S_1$ is 
incompatible with Proposition 2.4 (iii), and the cases $\rdescent{w} = S_0$ 
and $\rdescent{w} = S_1$ are incompatible with Proposition 2.4 (iv).
Because of the assumption about $\ldescent{w}$, the situation of Proposition 
2.4 (i) can only occur if $w = I_0$ or $w = I_1$, and these are special
cases of Proposition 2.4 (ii).  The result now follows by elimination.
\qed\enddemo

\head 3. Proof of main result \endhead

Our main result is the following

\proclaim{Theorem 3.1}
Let $W$ be a Coxeter group of type $\widetilde{A}_{n-1}$, let $w \in W$ be
arbitrary, and let $x \in W_c$.  Then $\mu(x, w) \in \{0, 1\}$.
\endproclaim

\demo{Proof}
Suppose first that $w$ is not fully commutative.

If $\ldescent{w}$ is not commutative, then Lemma 2.5 (i) shows that we
must have $s \in \ldescent{w} \backslash \ldescent{x}$ and the result follows
from Proposition 2.1 (iv).  Similarly, if $\rdescent{w}$ is not commutative,
then the result follows from Lemma 2.5 (i) and Proposition 2.1 (v).

Now let $w = w_0, w_1, \ldots, w_k$ be the sequence given in Proposition 2.3.  
We proceed by induction on $k$; the previous paragraph deals with the
case $k = 0$.  Let $I = \{s, t\}$ be a pair of noncommuting generators with 
the property that $w_1 = {^*w}_0$ with respect to $I$.  We may assume 
without loss of generality that $s \in \ldescent{w}$.  If $s \not\in 
\ldescent{x}$, then we are done by Proposition 2.1 (iv), so we may assume 
that $s \in \ldescent{x}$.  By Lemma 2.5 (i), we have $t \not\in 
\ldescent{x}$.  Since $\ldescent{x} \cap I$ consists of a single element, 
we can apply a left star operation to $x$ with respect to $I$.  By 
Proposition 2.1 (ii), we have $\mu(x, w) = \mu(^*x, {^*w})$.  We now repeat
the argument with $^*w$ in place of $w$, and the result follows by
induction.

We may now suppose that $w \in W_c$.  Proposition 2.4 shows that there
are four cases to consider.  

The first case is that $w$ is a product of commuting generators.  In order 
for $\mu(x, w)$ not to be zero, we need $x < w$, which means that $x$ is a 
product of a proper subset of the aforementioned commuting generators, and 
that there exists $s \in \ldescent{w} \backslash \ldescent{x}$.  The result 
then follows from Proposition 2.1 (iv).

The second case is that $n$ is even and $w$ is an alternating product of 
$I_0$ and $I_1$.  We may assume that $\ldescent{w} \subseteq \ldescent{x}$ 
and $\rdescent{w} \subseteq \rdescent{x}$, or we are done by Proposition 
2.1 (iv) and (v).  Since $\ldescent{x}$ and $\rdescent{x}$ are commutative,
we must have $\ldescent{w} = \ldescent{x}$ 
(because $\ldescent{w} \in \{S_0, S_1\}$) 
and $\rdescent{w} = \rdescent{x}$
(because $\rdescent{w} \in \{S_0, S_1\}$).
By Lemma 2.5 (ii), $w$ and $x$ must both be alternating products of $I_0$ and
$I_1$, and furthermore, both these alternating products must (a) start with
the same $I_i$, and (b) end with the same $I_j$.  It follows that
$\ell(w) = \ell(x) \mod 2$, which proves that $\mu(x, w) = 0$ by Proposition
2.1 (i), and the result follows.

The third case is the situation of Proposition 2.4 (iii).  We may assume
that $s \in \ldescent{x}$, or we are done by Proposition 2.1 (iv).  We also
have $t \not\in \ldescent{x}$ by Lemma 2.5 (i).  Since $\ldescent{x} \cap I$ 
is a singleton, the element $^*x$ with respect to $I$ is defined.  By
Proposition 2.1 (ii), we have $\mu(x, w) = \mu(^*x, {^*w})$.  By the 
assumptions on Proposition 2.4 (iii), $\ell({^*w}) = \ell(w) - 1$, and the 
proof is completed by induction on $\ell(w)$.

The fourth case is the situation of Proposition 2.4 (iv), and we argue
as in the previous paragraph, using Proposition 2.1 (iii) and (v).
\qed\enddemo

\remark{Remark 3.2}
It follows easily from Proposition 2.1 (i), (iv), (v) that Theorem 3.1 is 
also true in the case $n = 2$.
\endremark

\head Acknowledgement \endhead

I thank Brant Jones for some helpful conversations and correspondence.

\leftheadtext{} \rightheadtext{}

\end